# UNIVERSAL OPTIMALITY OF PATTERSON'S CROSSOVER DESIGNS[1]

By Kirti R. Shah, Mausumi Bose and Damaraju Raghavarao

*University of Waterloo, Indian Statistical Institute and Temple University*

We show that the balanced crossover designs given by Patterson [*Biometrika* **39** (1952) 32–48] are (a) universally optimal (UO) for the joint estimation of direct and residual effects when the competing class is the class of connected binary designs and (b) UO for the estimation of direct (residual) effects when the competing class of designs is the class of connected designs (which includes the connected binary designs) in which no treatment is given to the same subject in consecutive periods. In both results, the formulation of UO is as given by Shah and Sinha [Unpublished manuscript (2002)].

Further, we introduce a functional of practical interest, involving both direct and residual effects, and establish (c) optimality of Patterson's designs with respect to this functional when the class of competing designs is as in (b) above.

**1. Introduction.** Crossover designs (repeated measurement designs or change-over designs) in $v$ treatments on $n$ experimental units in $p$ periods are useful in a broad spectrum of research areas, including agriculture [2], dairy husbandry [3], bioassay procedures [4], clinical trials [5], psychological experiments [8] and weather modification experiments [17]. The advantages of the crossover design are its cost and the elimination of interunit variability. In the following, we assume that each treatment produces a direct effect in the period of its application and a residual effect in the subsequent period of its application.

Williams [23] gave designs for $p = v$ which were balanced in the sense that every paired difference of direct (residual) effects was estimated with the same precision. Patterson [18] gave combinatorial conditions for balance and

Received September 2002; revised February 2005.
[1]Supported by the Natural Sciences and Research Council of Canada in the form of a research grant.

also gave a number of methods for construction of such designs when $p \leq v$ and when $n$ is as small as possible. Since $p$ and $n$ are small, these designs are very attractive to practitioners. All these designs had the property that no treatment immediately succeeds itself on the same subject.

Hedayat and Afsarinijad [6] showed that when $p = v$, a balanced design is universally optimal (UO) (as defined in [9]) for estimation of the direct (residual) effects when the designs in the competing class are uniform on periods as well as subjects. Cheng and Wu [1] showed that these designs are UO for the estimation of residual effects when the competing designs may not be uniform over subjects or periods, but again no treatment succeeds itself on the same subject. Kunert [10] showed that when $n = vt$, a balanced uniform design is UO for direct effects if $v \geq 3$ and $t = 1$ or if $v \geq 6$ and $t = 2$. Hedayat and Yang [7] generalized this to the case where $v \geq 3$ and $t \leq (v-1)/2$. The results of Kunert [10] and of Hedayat and Yang [7] were proved without any condition on the competing designs. However, there do not appear to be any available results on the optimality of balanced crossover designs when $p < v$.

Cheng and Wu [1] also introduced what are called strongly balanced designs where each of the $v^2$ pairs of treatments occurs in consecutive periods for the same subject an equal number of times. They established some strong optimality properties for these designs. However, these designs require $p = vt$ or $vt + 1$ and also require $n$ to be large.

Kushner [14] gave a novel approximate design theory approach to obtain UO designs for arbitrary values of $p$ and $v$. Further, Kushner [15] gave exact designs which are UO for direct effects for every pair $(v, p)$ for some $n$.

Kushner's results are very attractive because they do not put any conditions on the competing designs. Their main limitation is that the values of $p$ or of $n$ are large. Further, in almost all cases these optimal designs are nonbinary (on the subjects). An attractive property of the binary balanced designs is that they are optimal when the residual effects are negligible [9].

Some authors have obtained optimal designs under different models. Kunert and Martin [11] gave optimal designs under an interference model. Kunert and Martin [12] considered models with correlated errors.

Kunert and Stufken [13] introduced a model where the residual effect of a treatment on itself is different from the residual effect when the treatment is followed by another treatment. An excellent review of the literature in this broad area up to 1996 is given by Stufken [22].

The balanced designs given by Williams [23] and by Patterson [18] are very attractive because they have a small number of periods and often involve a small to moderate number of subjects. These designs have been around for a long time and are generally believed to be efficient. However, precise optimality results are rather limited in nature.

In this paper we establish some strong optimality properties of the Patterson designs. We first show that, within the class of binary designs, these designs are UO (in the sense of Shah and Sinha [21]) for the joint estimation of direct and residual effects. This is a very strong property because it implies UO for the estimation of the direct (residual) effects and a great deal more. For the rest of the paper we refer to UO as formulated in [21].

Next, we establish the UO property of these designs for the estimation of direct (residual) effects where the only restriction on a competing design is that no treatment immediately succeeds itself on the same subject. We also compute lower bounds for the efficiencies of these designs within the *unrestricted* class of competing designs and find that these are very high (0.99 or higher), giving rise to the speculation that, when the fully efficient designs as described by Kushner [15] do not exist, these designs might, in fact, be optimal for specific criteria, such as A-optimality.

Further, we introduce an optimality function of practical interest and show that a Patterson crossover design is optimal for this functional, again with the restrictions on the competing designs that they are connected and that no treatment immediately succeeds itself on the same subject.

In this paper we compare the information matrix for a Patterson design with the average (over permutations of treatment labels) of the information matrix for a competing design. This average has a form which is much simpler than the original matrix. Further, we first prove the optimality results for the model without period effects and then for the model with period effects. This works well because, when we introduce period effects, the information matrix for a Patterson design is unaltered, whereas, for a competing design, it is reduced by a nonnegative definite matrix.

**2. Preliminaries.** Let us consider crossover designs where $v$ treatments are arranged in $p$ rows and $n$ columns. The rows correspond to periods whereas the columns correspond to the subjects. A crossover design is said to be balanced if we have the following:

(a) It is uniform over periods, that is, every treatment occurs $t$ times in each period.

(b) The design with subjects as blocks forms a balanced incomplete block design (BIBD).

(c) The design with subjects as blocks and last period omitted also forms a BIBD.

(d) Every ordered pair of distinct treatments occurs in consecutive periods in units $\lambda$ times (any pair of identical treatments does not occur in consecutive periods).

(e) In the set of $t$ subjects receiving a particular treatment in the last period, every other treatment is applied $\lambda$ times in the first $(p-1)$ periods.

These conditions are equivalent to conditions I–VII given by Patterson. We shall call such designs Patterson designs. Existence of a Patterson design implies

$$p \leq v, \qquad n = vt, \qquad \lambda = t(p-1)/(v-1),$$

where $t$ and $\lambda$ are positive integers. We shall assume that the parameter values $p, v, n$ are such that a Patterson design exists.

The crossover designs given by Williams [23] are balanced with $p = v$. However, one can often find balanced designs with $p < v$ and where $n$ is not too large. For $v = 4$ and $p = 3$, Patterson ([18], Figure 3) gives the following design with $n = 12$:

```
a  b  c  b  d  a  d  a  c  c  d  b
b  c  a  a  b  d  a  c  d  b  c  d
c  a  b  d  a  b  c  d  a  d  b  c
```

Patterson [18] gave several methods of construction for balanced crossover designs with $p \leq v$. Our Table 1, above extracted from Table 1 in [18] giving designs available for small values of $p$, is of practical interest.

Table 1 gives designs with the minimum value of $n$ for given values of $p$ and $v$ when $n \leq 60$. Many more designs can be constructed using the various methods given by Patterson. Further, such designs are available with $p = v = n$ when $v$ is even and $p = v = n/2$ when $v$ is odd [23]. Thus, this is a rich class of designs and it contains many designs of interest to practitioners.

We exclude the case $p = v = 2, n = 2t$ even though for these parameter values designs satisfying the combinatorial conditions exist. This is because in this case neither the direct nor the residual effects are estimable.

Let $d(i, j)$ be the treatment assigned to the $j$th subject in the $i$th period and let $y_{ij}$ denote the response obtained from that subject in that period. We assume that the $y_{ij}$'s are uncorrelated with common variance $\sigma^2$ and

$$(2.1) \qquad E(y_{ij}) = \mu + \alpha_i + \gamma_j + \tau_{d(i,j)} + \delta_{d(i-1,j)},$$

$i = 1, 2, \ldots, p; j = 1, 2, \ldots, n; \delta_{d(0,j)} = 0$ for all $j$, where $E(\cdot)$ denotes the expected value of the variable in the parentheses, $\mu$ is the general mean, and $\alpha, \gamma, \tau$ and $\delta$ are the period, subject, direct and residual treatment effects, respectively.

A crossover design is said to be connected if $\tau_i - \tau_{i'}$ and $\delta_i - \delta_{i'}$ are estimable for $i \neq i'$. All of Patterson's designs are connected. Let $\mathcal{D}$ denote the class of connected crossover designs using $n \ (= vt)$ subjects for comparing $v$ treatments in $p$ periods, with the restriction that, in each column of a design, adjoining positions are occupied by distinct treatments. Further, let $\mathcal{B}$ denote the subclass of $\mathcal{D}$ consisting of designs which are binary in the sense that no treatment is applied more than once to any subject.

We define vectors $\boldsymbol{\alpha}(p \times 1), \boldsymbol{\gamma}(n \times 1), \boldsymbol{\tau}(v \times 1)$ and $\boldsymbol{\delta}(v \times 1)$ whose components represent the above effects:

Table 1

| $p$ | 3 | 3 | 3 | 3 | 4 | 4 | 4 | 4 | 4 | 5 | 5 | 5 | 5 | 5 | 6 | 6 | 6 | 6 |
|---|---|---|---|---|---|---|---|---|---|---|---|---|---|---|---|---|---|---|
| $v$ | 3 | 7 | 8 | 11 | 4 | 5 | 7 | 8 | 13 | 5 | 7 | 8 | 11 | 13 | 6 | 7 | 8 | 11 |
| $n$ | 6 | 21 | 56 | 55 | 4 | 20 | 14 | 56 | 52 | 10 | 21 | 56 | 55 | 39 | 6 | 42 | 56 | 22 |

Let

$n_{iu} =$ number of appearances of treatment $i$ on subject $u$,

$\tilde{n}_{iu} =$ number of appearances of treatment $i$ on subject $u$
in the first $(p-1)$ periods,

$s_{ij} =$ number of appearances of treatment $i$ preceded
by treatment $j$ on the same unit,

$l_{ik} =$ number of appearances of treatment $i$ in period $k$,

$\tilde{l}_{i1} = 0, \qquad \tilde{l}_{ik} = l_{i(k-1)} \qquad \text{for } k \geq 2.$

We now define the frequency matrices $\mathbf{N} = (n_{iu})$, $\tilde{\mathbf{N}} = (\tilde{n}_{iu})$, $\mathbf{S} = (s_{ij})$, $\mathbf{L} = (l_{ik})$, $\tilde{\mathbf{L}} = (\tilde{l}_{ik})$. Further, let $\text{diag}(\mathbf{r})$ denote the $v \times v$ diagonal matrix whose elements are the replication numbers for the $v$ treatments in the entire design. Also, let $\text{diag}(\tilde{\mathbf{r}})$ denote the $v \times v$ diagonal matrix whose elements are the replication numbers for the treatments in the first $(p-1)$ periods only.

The information matrix for $(\boldsymbol{\tau}, \boldsymbol{\delta}, \boldsymbol{\alpha}, \boldsymbol{\gamma})$ is given by

$$(2.2) \qquad \mathcal{I}(\boldsymbol{\tau}, \boldsymbol{\delta}, \boldsymbol{\alpha}, \boldsymbol{\gamma}) = \begin{pmatrix} \text{diag}(\mathbf{r}) & \mathbf{S} & \mathbf{L} & \mathbf{N} \\ \cdot & \text{diag}(\tilde{\mathbf{r}}) & \tilde{\mathbf{L}} & \tilde{\mathbf{N}} \\ \cdot & \cdot & n\mathbf{I}_p & \mathbf{J}_{pn} \\ \cdot & \cdot & \cdot & p\mathbf{I}_n \end{pmatrix},$$

where $\mathbf{I}_a$ denotes the identity matrix of order $a$ and $\mathbf{J}_{ab}$ denotes an $a \times b$ matrix with all elements unity. (See equation (2.5) in [1].)

The information matrix for $(\boldsymbol{\tau}, \boldsymbol{\delta}, \boldsymbol{\alpha})$ eliminating $\boldsymbol{\gamma}$ is given by

$$(2.3)\ \mathcal{I}(\boldsymbol{\tau}, \boldsymbol{\delta}, \boldsymbol{\alpha}|\boldsymbol{\gamma}) = \begin{pmatrix} \text{diag}(\mathbf{r}) - \frac{1}{p}\mathbf{N}\mathbf{N}^t & \mathbf{S} - \frac{1}{p}\mathbf{N}\tilde{\mathbf{N}}^t & \mathbf{L} - \frac{1}{p}\mathbf{N}\mathbf{J}_{np} \\ \cdot & \text{diag}(\tilde{\mathbf{r}}) - \frac{1}{p}\tilde{\mathbf{N}}\tilde{\mathbf{N}}^t & \tilde{\mathbf{L}} - \frac{1}{p}\tilde{\mathbf{N}}\mathbf{J}_{np} \\ \cdot & \cdot & n\mathbf{I}_p - \frac{n}{p}\mathbf{J}_p \end{pmatrix}.$$

Here $\mathbf{J}_k$ is a $k \times k$ matrix with all elements unity.

When the period effects are ignored, the information matrix for $(\boldsymbol{\tau}, \boldsymbol{\delta})$ eliminating $\boldsymbol{\gamma}$ is seen to be

$$\mathcal{I}(\boldsymbol{\tau}, \boldsymbol{\delta}|\boldsymbol{\gamma}) = \begin{pmatrix} \operatorname{diag}(\mathbf{r}) - \dfrac{1}{p}\mathbf{NN}^t & \mathbf{S} - \dfrac{1}{p}\mathbf{N}\tilde{\mathbf{N}}^t \\ \cdot & \operatorname{diag}(\tilde{\mathbf{r}}) - \dfrac{1}{p}\tilde{\mathbf{N}}\tilde{\mathbf{N}}^t \end{pmatrix}. \tag{2.4}$$

Using $\begin{pmatrix} \frac{1}{n}\mathbf{I}_p & 0 \\ -\frac{1}{np}\mathbf{J}_{np} & \frac{1}{p}\mathbf{I}_n \end{pmatrix}$ as a $g$-inverse of $\begin{pmatrix} n\mathbf{I}_p & \mathbf{J}_{pn} \\ \mathbf{J}_{np} & p\mathbf{I}_n \end{pmatrix}$, we get the information matrix for $(\boldsymbol{\tau}, \boldsymbol{\delta})$ eliminating $(\boldsymbol{\alpha}, \boldsymbol{\gamma})$ as

$$\begin{aligned}\mathcal{I}(\boldsymbol{\tau}, \boldsymbol{\delta}|\boldsymbol{\alpha}, \boldsymbol{\gamma}) \\ = \mathcal{I}(\boldsymbol{\tau}, \boldsymbol{\delta}|\boldsymbol{\gamma}) - \begin{pmatrix} \dfrac{1}{n}\mathbf{LL}^t - \dfrac{1}{np}\mathbf{NJ}_{np}\mathbf{L}^t & \dfrac{1}{n}\mathbf{L}\tilde{\mathbf{L}}^t - \dfrac{1}{np}\mathbf{NJ}_{np}\tilde{\mathbf{L}}^t \\ \cdot & \dfrac{1}{n}\tilde{\mathbf{L}}\tilde{\mathbf{L}}^t - \dfrac{1}{np}\tilde{\mathbf{N}}\mathbf{J}_{np}\tilde{\mathbf{L}}^t \end{pmatrix}.\end{aligned} \tag{2.5}$$

For a Patterson design, $\mathbf{L} = (t\mathbf{J}_{vp})$, $\tilde{\mathbf{L}} = (\mathbf{0}|t\mathbf{J}_{v,p-1})$, $\operatorname{diag}(\mathbf{r}) = pt\mathbf{I}_v$, $\operatorname{diag}(\tilde{\mathbf{r}}) = t(p-1)\mathbf{I}_v$. Further, $\mathbf{NN}^t = p(t-\lambda)\mathbf{I}_v + p\lambda\mathbf{J}_v$, $\tilde{\mathbf{N}}\tilde{\mathbf{N}}^t = ((p-1)t - (p-2)\lambda)\mathbf{I}_v + (p-2)\lambda\mathbf{J}_v$ and $\mathbf{S} = \lambda(\mathbf{J}_v - \mathbf{I}_v)$. We note that, for all designs in the design class $\mathcal{D}$, the diagonal elements of $\mathbf{S}$ are all zeros.

Without loss of generality, we arrange the subjects so that the first $n_1$ subjects have treatment 1 in the last period, the next $n_2$ units have treatment 2 in the last period and so on.

This permits us to see that

$$\tilde{\mathbf{N}}^t = \mathbf{N}^t - \begin{bmatrix} \mathbf{1}_{n_1} & \mathbf{0} & \cdots & \mathbf{0} \\ \mathbf{0} & \mathbf{1}_{n_2} & \cdots & \mathbf{0} \\ \cdots & \cdots & \cdots & \cdots \\ \mathbf{0} & \mathbf{0} & \cdots & \mathbf{1}_{n_v} \end{bmatrix},$$

where $\mathbf{1}_h$ denotes an $h \times 1$ matrix with all elements unity.

This gives us

$$\begin{aligned}\mathbf{N}\tilde{\mathbf{N}}^t &= \mathbf{NN}^t - \boldsymbol{\Theta} \quad \text{and} \\ \tilde{\mathbf{N}}\tilde{\mathbf{N}}^t &= \mathbf{NN}^t - \boldsymbol{\Theta} - \boldsymbol{\Theta}^t + \operatorname{\mathbf{diag}}(n_1, n_2, \ldots, n_v),\end{aligned} \tag{2.6}$$

where $\boldsymbol{\Theta} = [\boldsymbol{\theta}_1, \boldsymbol{\theta}_2, \ldots, \boldsymbol{\theta}_v]$. Here, $\boldsymbol{\theta}_i$ is the sum of the $n_i$ columns of $\mathbf{N}$ corresponding to the $n_i$ subjects where treatment $i$ is in the last period.

It is easy to verify that, for a Patterson design,

$$\boldsymbol{\Theta} = (t-\lambda)\mathbf{I}_v + \lambda\mathbf{J}_v, \qquad \mathbf{NN}^t = p\boldsymbol{\Theta},$$
$$\mathbf{N}\tilde{\mathbf{N}}^t = (p-1)\boldsymbol{\Theta} \quad \text{and} \quad \tilde{\mathbf{N}}\tilde{\mathbf{N}}^t = (p-2)\boldsymbol{\Theta} + t\mathbf{I}_v.$$

For a design $d$, the information matrix for direct effects, residual effects and period effects eliminating the subject effects, that is, $\mathcal{I}(\boldsymbol{\tau}, \boldsymbol{\delta}, \boldsymbol{\alpha}|\boldsymbol{\gamma})$

of (2.3), will be denoted by $\mathbf{M}_d$. Let $d^*$ denote a Patterson design. The matrix $\mathbf{M}_d$ for a Patterson design is given by

$$(2.7) \quad \mathbf{M}_{d^*} = \begin{bmatrix} \frac{vt(p-1)}{(v-1)}\mathbf{H} & -\frac{vt(p-1)}{p(v-1)}\mathbf{H} & 0 \\ -\frac{vt(p-1)}{p(v-1)}\mathbf{H} & \frac{t(p-1)(pv-v-1)}{p(v-1)}\mathbf{H} + \frac{t(p-1)}{pv}\mathbf{J}_v & \frac{t}{p}\mathbf{1}_v(-(p-1),\mathbf{1}_{p-1}^t) \\ 0 & \frac{t}{p}\begin{pmatrix}-(p-1)\\ \mathbf{1}_{p-1}\end{pmatrix}\mathbf{1}_v^t & n\mathbf{I}_p - \frac{n}{p}\mathbf{J}_p \end{bmatrix}.$$

Here $\mathbf{H} = \mathbf{I}_v - \mathbf{J}_v/v$.

Similarly, for $d \in \mathcal{D}$, let $\mathbf{C}_d$ denote the $2v \times 2v$ information matrix for direct and residual effects eliminating the subject and period effects. We can write $\mathbf{C}_d$ as

$$(2.8) \quad \mathcal{I}(\boldsymbol{\tau}, \boldsymbol{\delta}|\boldsymbol{\alpha}, \boldsymbol{\gamma}) = \mathbf{C}_d = \begin{pmatrix} \mathbf{C}_{d11} & \mathbf{C}_{d12} \\ \mathbf{C}_{d21} & \mathbf{C}_{d22} \end{pmatrix},$$

where $\mathbf{C}_{d11}$ corresponds to the part for direct effects. $\mathbf{C}_{d12}, \mathbf{C}_{d21}$ and $\mathbf{C}_{d22}$ are described similarly. Note that $\mathbf{C}_{d21} = \mathbf{C}_{d12}^t$. Sometimes in the sequel we shall drop the suffix $d$.

For the Patterson design $d^*$, let $\mathbf{C}_{ij}^* = \mathbf{C}_{d^*ij}$. Using (2.4), (2.5) and (2.8), we have

$$(2.9) \quad \begin{aligned} \mathbf{C}_{11}^* &= \frac{vt(p-1)}{(v-1)}\mathbf{H}, \qquad \mathbf{C}_{12}^* = -\frac{vt(p-1)}{p(v-1)}\mathbf{H} \quad \text{and} \\ \mathbf{C}_{22}^* &= \frac{t(p-1)(pv-v-1)}{p(v-1)}\mathbf{H}. \end{aligned}$$

**3. UO for joint estimation in the design class $\mathcal{B}$.** In this section we shall show that $d^*$ is universally optimal (UO) for the joint estimation of direct and residual effects when the designs in the competing class are connected and are binary over subjects, that is, $n_{iu} = 0$ or $1$.

For formulations of UO one is referred to Kiefer [9], Shah and Sinha [20] and Shah and Sinha [21]. Here, we shall use the formulation of Shah and Sinha [21], which may be described as follows.

Let $\mathbf{C}_d$ denote a $v \times v$ direct effects information matrix (resp., $v \times v$ residual effects information matrix; or $2v \times 2v$ joint direct-residual effects information matrix) of design $d$. Let $g$ be a permutation of $\{1, 2, \ldots, v\}$, that is, $g \in S_v$, the symmetric group on $\{1, 2, \ldots, v\}$. A design $d_0$ having information matrix $\mathbf{C}_{d_0}$ is said to be UO in an appropriate design class if it minimizes every real valued function $\phi(\mathbf{C})$ (defined on the set of nonnegative definite matrices) that satisfies the following conditions:

(1) $\phi(\mathbf{C}_{dg}) = \phi(\mathbf{C}_d)$, where $dg$ is the design obtained by permuting treatment labels according to $g$.
(2) $\mathbf{C}_d \geq \mathbf{C}_f \Rightarrow \phi(\mathbf{C}_d) \leq \phi(\mathbf{C}_f)$, where $d$ and $f$ are any two designs.

(3) $\phi(\sum w_g \mathbf{C}_{dg}) \leq \phi(\mathbf{C}_d)$, where $w_g$ are all rational weights satisfying $\sum w_g = 1$. Here $g$ runs over all permutations in $S_v$.

It may be noted that every convex functional satisfies (3). This formulation of UO is an extension of Kiefer's original formulation [9], in the sense that the condition of convexity is replaced by the slightly weaker condition (3) above. See [20] and also [21] for a discussion of this.

A sufficient condition due to Shah and Sinha [21] for $d_0$ to be UO is

(3.1) $$\sum w_g \mathbf{C}_{dg} \leq \mathbf{C}_{d_0} \qquad \text{for every } d.$$

In (3.1) the $w_g$'s can be any specific set of weights (which may depend upon $d$). In the sequel we will use (3.1) when $w_g = 1/v!$ for all $g \in S_v$.

Let $\mathbf{M}_{dg}$ denote the matrix obtained from $\mathbf{M}_d$ by permuting treatment labels according to $g$. We shall first show that

(3.2) $$\mathbf{M}_{d^*} = \sum_g \mathbf{M}_{dg}/v!.$$

To show this, we shall state the following lemma which is easily established.

LEMMA 3.1. *Let $\mathbf{A}$ be a $k \times k$ matrix and let $g$ denote a permutation on $\{1, 2, \ldots, k\}$ for which the permutation matrix is $\mathbf{P}_g$. Then $\bar{\mathbf{A}} = \sum_{g \in S_k} \mathbf{P}_g^t \mathbf{A} \mathbf{P}_g / k!$ is a completely symmetric matrix with diagonal and off-diagonal elements $a$ and $b$, respectively, given by*

$$a = \sum_i a_{ii}/k \quad \text{and} \quad b = \left(s - \sum_i a_{ii}\right) \Big/ k(k-1).$$

*Here $s = \sum_i \sum_j a_{ij}$ is the sum of all elements of $\mathbf{A}$.*

We now consider the various submatrices of $\mathbf{M}_d$ for a binary design and show that, for each of these, the average over the permutations of treatment labels equals the corresponding expression for $d^*$.

For any binary design, the $i$th diagonal element of $\mathbf{NN}^t$ is $\sum_i \sum_u n_{iu}^2 = \sum_i \sum_u n_{iu} = r_i$, the replication number for the $i$th treatment. (This does not hold for a nonbinary design.) The average of $r_i$ over all permutations is $pt$, which is the replication number for $d^*$. Further, for a binary design the $i$th diagonal element of $\mathbf{\Theta}$ is $n_i$. The average of $n_i$ over all permutations is $t = n/v$, the common diagonal element of $\mathbf{\Theta}$ for $d^*$.

Let $g$ be a permutation on $\{1, 2, \ldots, v\}$. We note that the matrix $\mathbf{M}_{dg}$ is given by $\mathbf{M}_{dg} = \mathbf{Q}_g^t \mathbf{M}_d \mathbf{Q}_g$, where $\mathbf{Q}_g = \begin{pmatrix} \mathbf{P}_g & 0 & 0 \\ 0 & \mathbf{P}_g & 0 \\ 0 & 0 & \mathbf{I} \end{pmatrix}$.

Using Lemma 3.1 and the expressions given in (2.3), one can easily verify that, if the design is binary and is connected for each of $\text{diag}(\mathbf{r}) - \mathbf{NN}^t/p, \mathbf{S} -$ $\mathbf{N}\tilde{\mathbf{N}}^t/p$ and for $\mathrm{diag}(\tilde{\mathbf{r}}) - \tilde{\mathbf{N}}\tilde{\mathbf{N}}^t/p$, the average over all $g$ is indeed the corresponding expression for $d^*$. We note here that, for the first two of these, the sum of all the elements is zero, whereas, for the third one, it is $vt(p-1)/p$.

The assumption of connectedness is crucial here. If the design is not connected, $\mathbf{N}\mathbf{N}^t$ and/or $\tilde{\mathbf{N}}\tilde{\mathbf{N}}^t$ would have a block diagonal form where the off-diagonal submatrices consist of zeros.

We also note that the average for $\mathbf{L}$ is $t\mathbf{J}_{vp}$ and the average for $\mathbf{N}$ is $(p/v)\mathbf{J}_{vn}$. Thus, the average for $\mathbf{L} - \mathbf{N}\mathbf{J}_{np}/p$ is $\mathbf{0}$.

Finally, we note that for each of $\tilde{\mathbf{L}}$ and $\tilde{\mathbf{N}}$ the average over the $v!$ permutations gives the corresponding expressions for a Patterson design. This completes the proof of the assertion.

Now we note that the adjustment for the period effects $\boldsymbol{\alpha}$ is equivalent to computing the Schur complement. Thus, the Schur complement of $\mathbf{M}_{d^*}$ is $\mathbf{C}_{d^*}$, whereas the Schur complement of $\mathbf{M}_{dg}$ is $\mathbf{C}_{dg}$. Here $\mathbf{C}_{d^*}$ and $\mathbf{C}_{dg}$ refer to the $2v \times 2v$ joint direct-residual effects information matrix. Since $\mathbf{M}_{d^*} = \sum_g \mathbf{M}_{dg}/v!$ and since the Schur complement is a concave function [19], we get

$$\mathbf{C}_{d^*} \geq \sum_g \mathbf{C}_{dg}/v!.$$

Using the sufficient condition (3.1), with weights $w_g = 1/v!, g \in S_v$, we see that $d^*$ is UO w.r.t. any design $d \in \mathcal{B}$.

As shown in [16] and in [21], UO for the joint estimation of two sets of parameters is a very strong property. In particular, it implies UO for the estimation of each set of parameters. It is shown in [21] that the converse is not true. A design can be UO for the estimation of the direct effects, as well as for the residual effects. However, it might fail to be UO for the joint estimation of the two.

**4. UO for direct (residual) effects in the design class $\mathcal{D}$.** We now consider the case where the competing class of designs is $\mathcal{D}$, a class of designs that contains all binary designs. Initially we shall assume that there are no period effects. We shall relax this assumption subsequently.

The $\mathbf{C}_d$ matrix for this case is given by the submatrices of the matrix $\mathbf{M}_d$ which correspond to the direct and the residual effects. These are the components of the information matrix $\mathbf{C}_d$ for the direct and the residual effects *ignoring* the period effects. We write the expressions for these dropping the suffix $d$,

(4.1)
$$\begin{aligned}
\mathbf{C}_{11} &= \mathrm{diag}(\mathbf{r}) - \mathbf{N}\mathbf{N}^t/p, \\
\mathbf{C}_{12} &= \mathbf{S} - \mathbf{N}\tilde{\mathbf{N}}^t/p, \quad \mathbf{C}_{21} = \mathbf{C}_{12}^t, \\
\mathbf{C}_{22} &= \mathrm{diag}(\tilde{\mathbf{r}}) - \tilde{\mathbf{N}}\tilde{\mathbf{N}}^t/p.
\end{aligned}$$

We also write below the expressions for these submatrices for $d^*$, the Patterson design,

$$
\begin{aligned}
\mathbf{C}_{11}^* &= \mathbf{C}_{d^*11} = \frac{vt(p-1)}{(v-1)}\mathbf{H}, \\
(4.2) \quad \mathbf{C}_{12}^* &= \mathbf{C}_{d^*12} = -\frac{vt(p-1)}{p(v-1)}\mathbf{H}, \qquad \mathbf{C}_{21}^* = (\mathbf{C}_{12}^*)^t, \\
\mathbf{C}_{22}^* &= \mathbf{C}_{d^*22} = \frac{t(p-1)(pv-v-1)}{p(v-1)}\mathbf{H} + \frac{t(p-1)}{pv}\mathbf{J}_v.
\end{aligned}
$$

Let $\bar{\mathbf{C}}$ denote the average of $\mathbf{C}$ over all permutations of treatment labels. To describe the structure of $\bar{\mathbf{C}}$, we introduce some notation. We define the following:

$$
\begin{aligned}
(4.3) \quad \beta &= \sum_i \sum_u n_{iu}^2, \\
l &= \sum_i (\text{sum of } n_{iu}\text{'s for subjects with treatment } i \text{ in the last period}).
\end{aligned}
$$

Using the expressions for $\mathbf{N}\tilde{\mathbf{N}}^t$ and $\tilde{\mathbf{N}}\tilde{\mathbf{N}}^t$ in terms of $\mathbf{N}\mathbf{N}^t$ and $\boldsymbol{\Theta}$ given in (2.6) and using Lemma 3.1, one can verify that the matrices $\bar{\mathbf{C}}_{11}, \bar{\mathbf{C}}_{12}$ and $\bar{\mathbf{C}}_{22}$ have the following structure:

$\bar{\mathbf{C}}_{11}$: Diagonal element is $(p^2vt - \beta)/pv$,
off-diagonal element is $-(p^2vt - \beta)/pv(v-1)$.
$\bar{\mathbf{C}}_{12}$: Diagonal element is $-(\beta - l)/pv$,
off-diagonal element is $(\beta - l)/pv(v-1)$.
$\bar{\mathbf{C}}_{22}$: Diagonal element is $(vt(p^2 - p - 1) - (\beta - sl))/pv$,
off-diagonal element is $(\beta - 2l + pvt(2-p))/pv(v-1)$.

To illustrate the nature of the computations, we consider $\bar{\mathbf{C}}_{22}$. The sum of all the elements of $\text{diag}(\tilde{\mathbf{r}}) - \tilde{\mathbf{N}}\tilde{\mathbf{N}}^t/p$ is $vt(p-1)/p$. The average of the diagonal elements of $\tilde{\mathbf{N}}\tilde{\mathbf{N}}^t = \mathbf{N}\mathbf{N}^t - \boldsymbol{\Theta} - \boldsymbol{\Theta}^t + \text{diag}(n_1, \ldots, n_v)$ is $(\beta - 2l + vt)/v$, whereas the average of the diagonal elements of $\text{diag}(\tilde{\mathbf{r}})$ is $t(p-1)$. Use of Lemma 3.1 yields the expressions given above.

From these we deduce that

$$
\begin{aligned}
(4.4) \quad \bar{\mathbf{C}}_{11} &= \frac{(p^2vt - \beta)}{p(v-1)}\mathbf{H}, \qquad \bar{\mathbf{C}}_{12} = -\frac{\beta - l}{p(v-1)}\mathbf{H} \quad \text{and} \\
\bar{\mathbf{C}}_{22} &= \frac{pvt(p-1) - (\beta - 2l) - t(v+p-1)}{p(v-1)}\mathbf{H} + \frac{t(p-1)}{pv}\mathbf{J}_v.
\end{aligned}
$$

We shall now work with the information for direct (residual) effects adjusted for residual (direct) effects and shall use condition (3.1) to show the UO property of Patterson designs.

We show below that both $(\bar{\mathbf{C}}_{22})^{-1}$ and $(\mathbf{C}^*_{22})^{-1}$ exist. It is then easy to see that $\mathbf{C}_{11.22} = \mathbf{C}_{11} - \mathbf{C}_{12}\mathbf{C}_{22}^{-1}\mathbf{C}_{21}$ for $\mathbf{C}^*$ and for $\bar{\mathbf{C}}$ are given by

$$\mathbf{C}^*_{11.22} = \frac{vt(p-1)}{(v-1)}\left\{1 - \frac{v}{p(pv-v-1)}\right\}\mathbf{H} \quad \text{and}$$

(4.5) $\quad \bar{\mathbf{C}}_{11.22} = \frac{1}{p(v-1)}\left\{p^2tv - \beta \right.$

$$\left. - \frac{(\beta - l)^2}{[pvt(p-1) - (\beta - 2l) - t(v+p-1)]}\right\}\mathbf{H}.$$

We note that $\mathbf{C}_{11.22}$ is the information matrix for the direct effects eliminating the residual effects. Further, $\bar{\mathbf{C}}_{11.22}$ is $(\sum_g \mathbf{C}_g/v!)_{11.22}$.

We now show that each of $\mathbf{C}^*_{22}$ and $\bar{\mathbf{C}}_{22}$ is nonsingular. We first consider $\bar{\mathbf{C}}_{22}$. Since $t(p-1)/pv \neq 0$, $\bar{\mathbf{C}}_{22}$ is nonsingular iff the coefficient of $\mathbf{H}$ in the expression for $\bar{\mathbf{C}}_{22}$ is nonzero. If this coefficient is zero, rank $\bar{\mathbf{C}}_{22}$ is unity. We now recall that the $(\delta_i - \delta_{i'})$'s are all estimable when the period effects are *eliminated*. These continue to be estimable when the period effects are ignored. Now, estimability of all $(\delta_i - \delta_{i'})$'s implies that rank $\mathbf{C}_{22} \geq v - 1$. This also implies that the sum of $\mathbf{C}_{22}$ over all permutations of treatment labels has rank at least $(v-1)$. We thus have rank $\bar{\mathbf{C}}_{22} \geq v - 1$ and, hence, the coefficient of $\mathbf{H}$ in the above expression for $\bar{\mathbf{C}}_{22}$ must be nonzero. From the expression (4.2) for $\mathbf{C}^*_{22}$, it is clear that it is of full rank.

We shall now show that $\mathbf{C}^*_{11.22} - \bar{\mathbf{C}}_{11.22}$ is n.n.d. To see this, let $l_{ij}$ denote the value of $n_{iu}$ in the $j$th of the $n_i$ subjects with treatment $i$ in the last period. Similarly, let $\beta_{ij}$ denote the contribution to $\beta$ for that subject. We now note that $\beta = \sum_{i=1}^{v} \sum_{j=1}^{n_i} \beta_{ij}, l = \sum_{i=1}^{v} \sum_{j=1}^{n_i} l_{ij}$. We note that $l \geq vt$. Further, for the $j$th subject with treatment $i$ in the last period, one $n_{iu}$ is $l_{ij}$ and there are $(v-1)$ other $n_{iu}$'s which add up to $p - l_{ij}$. Hence,

$$\beta_{ij} \geq l_{ij}^2 + (p - l_{ij}) = l_{ij}(l_{ij} - 1) + p.$$

It follows that

$$\beta_{ij} - 2l_{ij} \geq l_{ij}^2 - 3l_{ij} + p = (l_{ij} - 1)(l_{ij} - 2) + p - 2,$$
$$\beta_{ij} - l_{ij} \geq l_{ij}^2 - 2l_{ij} + p = (l_{ij} - 1)^2 + p - 1.$$

Since $l_{ij}$ is 0, 1, 2 or greater than 2, we have

$$\beta_{ij} - 2l_{ij} \geq p - 2, \qquad \beta_{ij} - l_{ij} \geq p - 1, \qquad \beta_{ij} \geq p.$$

This gives $\beta \geq pvt, \beta - l \geq vt(p-1), \beta - 2l \geq vt(p-2)$. Using the above relations, it is easy to see that $\mathbf{C}^*_{11.22} \geq \bar{\mathbf{C}}_{11.22}$, that is, $\mathbf{C}^*_{11.22} - \bar{\mathbf{C}}_{11.22}$ is n.n.d.

An application of (3.1) shows that $d^*$ is UO (compared with any design in $\mathcal{D}$) if

$$\mathbf{C}^*_{11.22} - \sum_g (\mathbf{C}_g)_{11.22}/v! \text{ is n.n.d.}$$

12 K. R. SHAH, M. BOSE AND D. RAGHAVARAOHere $(\mathbf{C}_g)_{11.22}$ is obtained by permuting treatment labels in $\mathbf{C}$ by $g$ and then computing the Schur complement $(\mathbf{C}_g)_{11.22}$. In turns out that the result is the same as applying permutations $g$ to the rows and columns $\mathbf{C}_{11.22}$, that is, $(\mathbf{C}_g)_{11.22} = (\mathbf{C}_{11.22})_g$.

Since the Schur complement $\mathbf{C}_{11.22}$ is a concave function [19],

$$\left(\sum_g \mathbf{C}_g\right)_{11.22} \geq \sum_g (\mathbf{C}_g)_{11.22}.$$

Since $\bar{\mathbf{C}}_{11.22} = (\frac{1}{v!} \sum_g \mathbf{C}_g)_{11.22}$, we have

$$\mathbf{C}^*_{11.22} \geq \bar{\mathbf{C}}_{11.22} = \left(\sum_g \mathbf{C}_g\right)_{11.22} \Big/ v! \geq \sum_g (\mathbf{C}_g)_{11.22}/v!.$$

Thus, $d^*$ is UO in $\mathcal{D}$ for the model without period effects.

We shall now introduce period effects and consider estimation of the direct and the residual effects eliminating the period effects. When we adjust for periods, $\mathbf{C}$ gets reduced by an n.n.d. matrix and, hence,

$$\mathbf{C} \text{ (adjusted for periods)} \leq \mathbf{C} \text{ (ignoring periods)}.$$

This implies ([19], Section 3.13)

$$\mathbf{C} \text{ (adjusted for periods)}_{11.22} \leq \mathbf{C} \text{ (ignoring periods)}_{11.22}.$$

We now show that, for $d^*$,

$$\mathbf{C}^* \text{ (adjusted for periods)}_{11.22} = \mathbf{C}^* \text{ (ignoring periods)}_{11.22}.$$

Note that, for $d^*$, the expressions for $\mathbf{C}_d$ given by (2.9) and (4.2) differ only in $\mathbf{C}^*_{22}$. For these two cases, the Moore–Penrose inverses are $\mathbf{C}^{*+}_{22} = \frac{p(v-1)}{t(p-1)(pv-v-1)}\mathbf{H}$ and $\mathbf{C}^{*+}_{22} = \frac{p(v-1)}{t(p-1)(pv-v-1)}\mathbf{H} + \frac{p}{tv(p-1)}\mathbf{J}_v$, respectively. Since $\mathbf{C}_{11.22} = \mathbf{C}_{11} - \mathbf{C}_{12}\mathbf{C}^+_{22}\mathbf{C}_{21}$, and since $\mathbf{H}\mathbf{J}_v = \mathbf{0}$, the result follows.

We now have

$$\mathbf{C}^* \text{ (adjusted for periods)}_{11.22} = \mathbf{C}^* \text{ (ignoring periods)}_{11.22}$$
$$\geq \frac{1}{v!} \sum_g (\mathbf{C}_g \text{ (ignoring periods)})_{11.22}$$
$$\geq \frac{1}{v!} \sum_g (\mathbf{C}_g \text{ (adjusted for periods)})_{11.22}.$$

This establishes the UO property of $d^*$ (relative to any design in $\mathcal{D}$) for the estimation of direct effects.

An analogous argument also works for the estimation of residual effects. We outline the relevant important steps here. The information matrix for residual effects eliminating direct effects is now given by $\mathbf{C}_{22.11} = \mathbf{C}_{22} - \mathbf{C}_{21}\mathbf{C}_{11}^{+}\mathbf{C}_{12}$, where $\mathbf{C}_{11}^{+}$ is the Moore–Penrose inverse of $\mathbf{C}_{11}$. The expressions for $\mathbf{C}_{22.11}^{*}$ and $\bar{\mathbf{C}}_{22.11} = (\sum_g \mathbf{C}_g/v)_{22.11}$ are given by

$$\mathbf{C}_{22.11}^{*} = \frac{t(p-1)(pv-v-1)}{p(v-1)}\mathbf{H} + \frac{t(p-1)}{pv}\mathbf{J}_v - \frac{vt(p-1)}{p^2(v-1)}\mathbf{H},$$

$$\bar{\mathbf{C}}_{22.11} = \frac{pvt(p-1) - (\beta - 2l) - t(v+p-1)}{p(v-1)}\mathbf{H}$$
$$+ \frac{t(p-1)}{pv}\mathbf{J}_v - \frac{(\beta - l)^2}{p(v-1)(p^2vt - \beta)}\mathbf{H}.$$

Use of $\beta \geq pvt, \beta - l \geq vt(p-1)$ and $\beta - 2l \geq vt(p-2)$ yields $\mathbf{C}_{22.11}^{*} \geq \bar{\mathbf{C}}_{22.11}$. As in the case of direct effects, one can now show the UO for the estimation of residual effects.

It should be noted that $d^*$ is *not* UO in the class $\mathcal{D}$ for the joint estimation of direct and residual effects. This has been shown in [21].

It is not known if a Patterson design using the smallest number of subjects for given values of $p$ and $v$ is UO for the estimation of direct (residual) effects in the whole class of competing designs with fixed values of $p, v$ and $n$. Hence, we obtain a lower bound for the efficiency factor for the estimation of direct effects along the lines of [15].

First, we note that, for an approximate optimal design, the information matrix $\mathbf{C}_{11.22}$ given by [15] is

$$\tilde{\mathbf{C}}_{11.22} = \frac{vt}{v-1}\left\{p - 1 - \frac{1}{p} - \frac{1}{p(p-1)v}\right\}\mathbf{H}.$$

$\phi(\tilde{\mathbf{C}}_{11.22})$ can serve as a lower bound to the $\phi$ value of the information matrix for any competing design. Further, we note that both $\tilde{\mathbf{C}}_{11.22}$ and $\mathbf{C}_{11.22}^{*}$ given by (4.5) are multiples of the matrix $\mathbf{H}$. Hence, a lower bound to the efficiency of a Patterson design for the estimation of the direct effects is given by

$$e^* = \left(1 - \frac{v}{p(pv-v-1)}\right) \bigg/ \left(1 - \frac{pv-v+1}{pv(p-1)^2}\right).$$

On simplification this reduces to

$$e^* = A/(A+v),$$

where $A = v^2(p-1)^2(pv(p-1) - p - v)$.

In Table 2 we give the values of $e^*$ for the 18 designs given in Table 1 of Section 2.

As remarked earlier, these are very high. Thus, while our optimality results are restricted to the design class $\mathcal{D}$, we see that Patterson designs have very high efficiencies even in the unrestricted class. As remarked in [15], when $p \leq v$, no design can be optimal in the unrestricted class for each of direct and residual effects. Patterson designs are UO for each of direct and residual effects within the class $\mathcal{D}$ and are highly efficient even when we do not put any restrictions on the class of competing designs.

**5. A functional of interest.** In this section we establish an interesting optimality property of $d^*$ when $\mathcal{D}$ is the class of competing designs.

Sometimes, one may wish to look for the best combination of direct and residual effects. Thus, we would like to compare the yield when treatment $j$ is followed by treatment $i$ with the yield when treatment $j'$ is followed by treatment $i'$. This means that we wish to estimate all functions of the form $\tau_i + \delta_j - \tau_{i'} - \delta_{j'}$ as precisely as possible. This leads to the minimization of

$$(5.1) \qquad A = \sum \mathrm{Var}(\hat{\tau}_i + \hat{\delta}_j - \hat{\tau}_{i'} - \hat{\delta}_{j'}),$$

where the summation is over $i, i', j, j'$ such that $i \neq j, i \neq i', i' \neq j', j \neq j'$.

Since the variances of estimable parametric functions are functions of $\mathbf{C}$, we may write $A$ as $A(\mathbf{C})$.

Let $\mathbf{C}^+ = \begin{pmatrix} \mathbf{\Omega} & \mathbf{\Theta} \\ \mathbf{\Theta}^t & \mathbf{\Delta} \end{pmatrix}$ denote the Moore–Penrose inverse of $\mathbf{C}$. It is easy to see from (2.5) that each of $\mathbf{\Omega}, \mathbf{\Theta}$ or $\mathbf{\Delta}$ has zero row (column) sums. Expressing $A$ as

$$A = \sum \{\mathrm{Var}(\hat{\tau}_i - \hat{\tau}_{i'}) + \mathrm{Var}(\hat{\delta}_j - \hat{\delta}_{j'}) + 2 \mathrm{Cov}(\hat{\tau}_i - \hat{\tau}_{i'}, \hat{\delta}_j - \hat{\delta}_{j'})\}$$

and using the above property of $\mathbf{\Omega}, \mathbf{\Theta}$ and $\mathbf{\Delta}$, one can show that

$$(5.2) \quad A = 2v \left[ \{(v-1) + (v-2)^2\} \left\{ \sum_1^v \alpha_{ii} + \sum_1^v \beta_{ii} \right\} - 2(v-1) \sum_1^v \theta_{ii} \right] \sigma^2,$$

where $\mathbf{\Omega} = (\alpha_{ij}), \mathbf{\Delta} = (\beta_{ij})$ and $\mathbf{\Theta} = (\theta_{ij})$. We note that (5.2) may not hold if row (column) sums of each $\mathbf{\Omega}, \mathbf{\Theta}$ and $\mathbf{\Delta}$ are not zero.

Table 2
*Efficiency lower bounds for Patterson designs for direct effects*

| $p$ | $v$ | $e^*$ | $p$ | $v$ | $e^*$ | $p$ | $v$ | $e^*$ |
|---|---|---|---|---|---|---|---|---|
| 3 | 3 | 0.993103 | 4 | 7 | 0.999783 | 5 | 11 | 0.999972 |
| 3 | 7 | 0.998885 | 4 | 8 | 0.999835 | 5 | 13 | 0.999980 |
| 3 | 8 | 0.999156 | 4 | 13 | 0.999939 | 6 | 6 | 0.999960 |
| 3 | 11 | 0.999563 | 5 | 5 | 0.999853 | 6 | 7 | 0.999971 |
| 4 | 4 | 0.999306 | 5 | 7 | 0.999930 | 6 | 8 | 0.999978 |
| 4 | 5 | 0.999565 | 5 | 8 | 0.999947 | 6 | 11 | 0.999988 |

One can see that $A$ can be written as $A = A(\mathbf{C}) = (\sum_{l \in \mathcal{L}} \boldsymbol{l}^t \mathbf{C}^+ \boldsymbol{l})\sigma^2$, where $\mathcal{L}$ denotes the set of coefficient vectors for all contrasts included in (5.1).

We shall first show that $A$ viewed as a functional on $\mathbf{C}$ and denoted by $A(\mathbf{C})$ satisfies the three conditions on $\phi(\mathbf{C})$ given in the UO formulation in Section 3.

To see that $A$ is invariant under a permutation of treatment labels, we shall show that a permutation $g$ changes a coefficient vector $\mathbf{l}$ into a vector $g(\mathbf{l}) = \mathbf{l}'$, where $\mathbf{l}'$ can be seen to satisfy all the necessary constraints. Let $\mathbf{P}_g$ be the $v \times v$ matrix for a permutation $g \in S_v$. Let $\mathbf{F}_g = \begin{pmatrix} \mathbf{P}_g & \mathbf{0} \\ \mathbf{0} & \mathbf{P}_g \end{pmatrix}$. Then $\mathbf{C}_{dg} = \mathbf{F}_g^t \mathbf{C}_d \mathbf{F}_g$, giving $\mathbf{C}_{dg}^+ = \mathbf{F}_g^t \mathbf{C}_d^+ \mathbf{F}_g$. Thus, $\mathbf{l}^t \mathbf{C}_{dg}^+ \mathbf{l} = \mathbf{l}'^t \mathbf{C}_d^+ \mathbf{l}'$, where $\mathbf{l}' = \mathbf{F}_g \mathbf{l}$. As $\mathbf{l}$ varies over $\mathcal{L}$, $\mathbf{l}'$ also varies over $\mathcal{L}$. Thus, condition (1) holds. Next, if $\mathbf{C}_1 \geq \mathbf{C}_2, \mathbf{C}_1^+ \leq \mathbf{C}_2^+$ and, hence, $\mathbf{l}^t \mathbf{C}_1^+ \mathbf{l} \leq \mathbf{l}^t \mathbf{C}_2^+ \mathbf{l}$ for all $\mathbf{l} \in \mathcal{L}$. Thus, condition (2) also holds.

To show that condition (3) holds, we first note that $\mathbf{C}^+$ is a convex function of $\mathbf{C}$, that is, $(\sum w_g \mathbf{C}_g)^+ \leq \sum w_g \mathbf{C}_g^+$, where the $w_g$'s are rational weights satisfying $\sum w_g = 1$. Thus, we have

$$\begin{aligned}
A\Big(\sum w_g \mathbf{C}_g\Big) &= \sum_{\mathbf{l}} \mathbf{l}^t \Big(\sum w_g \mathbf{C}_g\Big)^+ \mathbf{l} \cdot \sigma^2 \\
&\leq \sum_{\mathbf{l}} \mathbf{l}^t \sum_g w_g \mathbf{C}_g^+ \mathbf{l} \cdot \sigma^2 \\
&\leq \sum_g w_g \sum_{\mathbf{l}} \mathbf{l}^t \mathbf{C}_g^+ \mathbf{l} \cdot \sigma^2 \\
&= \sum_g w_g A(\mathbf{C}_g) = \sum_g w_g A(\mathbf{C}) = A(\mathbf{C}).
\end{aligned}$$

This completes the verification. If we take $w_g = 1/v!$, we get $A(\bar{\mathbf{C}}) \leq A(\mathbf{C})$. Thus, to show that $A(\mathbf{C}_{d^*}) \leq A(\mathbf{C}_d)$, it is enough to show that $A(\mathbf{C}_{d^*}) \leq A(\bar{\mathbf{C}}_d)$.

Again, we initially assume that there are no period effects in the model. We first express $\mathbf{C}_{d^*}$ and $\bar{\mathbf{C}}_d$ as

(5.3)
$$\mathbf{C}_{d^*} = \begin{pmatrix} a^* \mathbf{H} & b^* \mathbf{H} \\ b^* \mathbf{H} & c^* \mathbf{H} + e \mathbf{J}_v \end{pmatrix},$$
$$\bar{\mathbf{C}}_d = \begin{pmatrix} \bar{a}_d \mathbf{H} & \bar{b}_d \mathbf{H} \\ \bar{b}_d \mathbf{H} & \bar{c}_d \mathbf{H} + e \mathbf{J}_v \end{pmatrix},$$

where $a^*, b^*, c^*, e, \bar{a}_d, \bar{b}_d$ and $\bar{c}_d$ obtained from (4.2) and (4.4) are

$$a^* = \frac{vt(p-1)}{v-1}, \qquad b^* = -\frac{vt(p-1)}{p(v-1)},$$
$$c^* = \frac{t(p-1)(pv-v-1)}{p(v-1)}, \qquad e = \frac{t(p-1)}{pv},$$

$$\bar{a}_d = \frac{(p^2vt - \beta)}{p(v-1)}, \qquad \bar{b}_d = -\frac{(\beta - l)}{p(v-1)},$$

$$\bar{c}_d = \frac{pvt(p-1) - (\beta - 2l) - t(v+p-1)}{p(v-1)}.$$

For simplicity of notation, we shall drop the subscript $d$ from $\bar{a}_d, \bar{b}_d$ and $\bar{c}_d$.

Since $\bar{\mathbf{C}}$ is an average of nonnegative definite (n.n.d.) matrices, it is n.n.d. Also, $\bar{\mathbf{C}}\binom{\mathbf{1}_v}{\mathbf{0}} = \mathbf{0}$ and $\bar{\mathbf{C}}\binom{\mathbf{0}}{\mathbf{1}_v} = ev\binom{\mathbf{0}}{\mathbf{1}_v}$. Other eigenvalues of $\bar{\mathbf{C}}$ are obtained as follows. Let $\mathbf{u}$ be a $(v \times 1)$ vector satisfying $\mathbf{Hu} = \mathbf{u}$. It is easy to verify that $\binom{\mathbf{u}}{\alpha \mathbf{u}}$ is an eigenvector of $\bar{\mathbf{C}}$ if $\alpha = \{(\bar{c} - \bar{a}) \pm \sqrt{(\bar{a}+\bar{c})^2 - 4(\bar{a}\bar{c} - \bar{b}^2)}\}/2$. The corresponding eigenvalue is $\{(\bar{a}+\bar{c}) \pm \sqrt{(\bar{a}+\bar{c})^2 - 4(\bar{a}\bar{c} - \bar{b}^2)}\}/2$. Since there are $(v-1)$ orthonormal choices for $\mathbf{u}$, each of $\{(\bar{a}+\bar{c}) \pm \sqrt{(\bar{a}+\bar{c})^2 - 4(\bar{a}\bar{c} - \bar{b}^2)}\}/2$ is an eigenvalue of $\bar{\mathbf{C}}$ of multiplicity $v-1$.

Since $\tau_i - \tau_{i'}$ and $\delta_j - \delta_{j'}$ are estimable in the model eliminating period effects, they are also estimable in the model ignoring period effects. Thus, rank $\bar{\mathbf{C}} \geq 2(v-1)$.

If $\bar{\Delta} = \bar{a}\bar{c} - \bar{b}^2 < 0$, $\bar{\mathbf{C}}$ has a negative eigenvalue. If $\bar{\Delta} = 0$, rank $\bar{\mathbf{C}} = v$, which is also a contradiction. Thus, $\bar{\Delta} > 0$. Similarly, $\Delta^* = a^*c^* - b^{*2} > 0$.

Direct calculations yield

(5.4)
$$\mathbf{C}^{*+} = \begin{pmatrix} c^*\mathbf{H}/\Delta^* & -b^*\mathbf{H}/\Delta^* \\ -b^*\mathbf{H}/\Delta^* & a^*\mathbf{H}/\Delta^* + \mathbf{J}_v/ev^2 \end{pmatrix},$$
$$\bar{\mathbf{C}}^+ = \begin{pmatrix} \bar{c}\mathbf{H}/\bar{\Delta} & -\bar{b}\mathbf{H}/\bar{\Delta} \\ -\bar{b}\mathbf{H}/\bar{\Delta} & \bar{a}\mathbf{H}/\bar{\Delta} + \mathbf{J}_v/ev^2 \end{pmatrix}.$$

Since the last $v$ components in each vector $\mathbf{l}$ are $v-2$ zeros, $+1$ and $-1$ (in some order), we can, in computing $\sum_{l \in \mathcal{L}} \mathbf{l}^t \mathbf{C}^+ \mathbf{l}$, ignore the term $\mathbf{J}_v/ev^2$ in $\mathbf{C}^{*+}$ and $\bar{\mathbf{C}}^+$. The $v \times v$ submatrices of the remaining matrix have zero row and column sums. Hence (5.2) is applicable to this matrix.

We can thus express $A^* = A(\mathbf{C}_{d^*})$ and $\bar{A} = A(\bar{\mathbf{C}}_d)$ as

$$A^* = 2v(v-1)\left[\{(v-1) + (v-2)^2\}\frac{(a^* + c^*)}{\Delta^*} + 2(v-1)\frac{b^*}{\Delta^*}\right]\sigma^2$$

and

$$\bar{A} = 2v(v-1)\left[\{(v-1) + (v-2)^2\}\frac{(\bar{a} + \bar{c})}{\bar{\Delta}} + 2(v-1)\frac{\bar{b}}{\bar{\Delta}}\right]\sigma^2.$$

If we write $\beta - 2l = vt(p-2) + x$ and $l = vt + y$, we can express $\bar{a}, \bar{b}$ and $\bar{c}$ as

$$\bar{a} = a^* - \frac{x + 2y}{p(v-1)}, \qquad \bar{b} = b^* - \frac{x+y}{p(v-1)}, \qquad \bar{c} = c^* - \frac{x}{p(v-1)}.$$

It was noted in Section 4 that $\beta - 2l \geq pvt(p-2)$ and $l \geq vt$. Hence, we have $x \geq 0$ and $y \geq 0$.

We may write $A^* = A(0,0)$ and $\bar{A} = A(x,y)$. To show that $A(0,0) \leq A(x,y)$, we shall proceed as follows:

Using the expressions for $\bar{a}, \bar{b}$ and $\bar{c}$ given above, it can be seen that

$$A(x,y) = \{(c_{11} + c_{12}x + c_{13}y)/(c_{21} + c_{22}x + c_{23}y - y^2/t(p-1))\}\sigma^2,$$

where

$$c_{11} = 2pv(v-1)(2pv^3 - 6pv^2 + 6pv - v^3 + 2v - 3),$$
$$c_{12} = c_{13} = -4pv(v-1)(v^2 - 2v + 2)/t(p-1),$$
$$c_{21} = tv(p-1)(p^2v - pv - p - v),$$
$$c_{22} = -(2pv + v - 1), \qquad c_{23} = -2(pv-1).$$

It is easy to verify that

$$\bar{\Delta} = \frac{t(p-1)}{p^2(v-1)^2}\{c_{21} + c_{22}x + c_{23}y - y^2/t(p-1)\}.$$

Since $\bar{\Delta} > 0$, it follows that $c_{21} + c_{22}x + c_{23}y - y^2/t(p-1) > 0$. Similarly, $\Delta^* > 0$ implies $c_{21} > 0$.

It follows that $A(x,y) - A(0,0)$ is strictly positive iff $(c_{21}c_{12} - c_{11}c_{22})x + (c_{21}c_{13} - c_{11}c_{23})y + c_{11}y^2/t(p-1)$ is strictly positive. For $p \geq 3, v \geq p$, the co-efficients of $x, y$ and $y^2$ are all seen to be strictly positive and, hence, $A(x,y) - A(0,0) > 0$ if $(x,y) \neq (0,0)$. When $p = 2$, we must have $y = 0$. We shall comment on this case later in this section.

We have thus shown that $A(x,y) \geq A(0,0)$ when period effects are ignored.

When we take period effects into the model, $\mathbf{C}_{d^*}$ in (5.3) gets reduced by $\begin{pmatrix} 0 & 0 \\ 0 & e\mathbf{H} \end{pmatrix}$ which has no effect on $A$. For $\bar{\mathbf{C}}$ we argue as follows. We first note that

$$\mathbf{C} \text{ (adjusted for periods)} \leq \mathbf{C} \text{ (ignoring periods)},$$

Since $\bar{\mathbf{C}}$ is obtained by averaging $\mathbf{C}_g$ over all permutations, it follows that

$$\bar{\mathbf{C}} \text{ (adjusted for periods)} \leq \bar{\mathbf{C}} \text{ (ignoring periods)},$$

$$\bar{\mathbf{C}}^+ \text{ (adjusted for periods)} \geq \bar{\mathbf{C}}^+ \text{ (ignoring periods)}.$$

Since $A(\mathbf{C}) = \sum_{\mathbf{l}} \mathbf{l}^t \mathbf{C}^+ \mathbf{1} \cdot \sigma^2$, it follows that adjustment for periods cannot decrease the value of $A(\mathbf{C})$.

We now summarize the situation as follows. Here, "adj" means adjusted for periods and "ign" means ignoring periods. We have seen that

$$A(\mathbf{C}_{d^*}(\text{adj})) = A(\mathbf{C}_{d^*}(\text{ign})),$$
$$A(\mathbf{C}_{d^*}(\text{ign})) \leq A(\bar{\mathbf{C}}_d(\text{ign})), \qquad d \in \mathcal{D},$$
$$A(\bar{\mathbf{C}}_d(\text{ign})) \leq A(\bar{\mathbf{C}}_d(\text{adj})).$$

These imply

$$A(\mathbf{C}_{d^*}(\mathrm{adj})) \le A(\bar{\mathbf{C}}_d(\mathrm{adj})), \qquad d \in \mathcal{D}.$$

This completes the proof of optimality of $d^*$ for the functional $A(\mathbf{C})$ in the design class $\mathcal{D}$.

In the definition of $A(\mathbf{C})$ we could also permit $i = i'$, as this would only add comparisons of the type $\delta_j - \delta_{j'}$. We have already seen that, for the estimation of residual effects, $d^*$ is UO in $\mathcal{D}$. Similarly, we could also permit $j = j'$.

It should be noted that the above proof was needed only when $d$ is nonbinary with $y > 0$. When the design is binary, or nonbinary with $y = 0$, the result follows from the UO property of $d^*$ (Section 3 of this paper; [21]) for the joint estimation of direct and residual effects.

**Acknowledgments.** The authors have greatly benefited from several discussions with Professor Bikas Kumar Sinha of the Indian Statistical Institute, Kolkata. His various suggestions are gratefully acknowledged. The authors are also grateful to an Associate Editor and to a referee for correcting an error and for making several suggestions which resulted in an improved presentation.

## REFERENCES

[1] Cheng, C.-S. and Wu, C.-F. (1980). Balanced repeated measurements designs. *Ann. Statist.* **8** 1272–1283. [Correction (1988) **11** 349.] MR0594644
[2] Cochran, W. G. (1939). Long-term agricultural experiments. *J. Roy. Statist. Soc. Suppl.* **6** 104–148.
[3] Cochran, W. G., Autrey, K. M. and Canon, C. Y. (1941). A double change-over design for dairy cattle feeding experiments. *J. Dairy Sci.* **24** 937–951.
[4] Finney, D. J. (1956). Cross-over designs in bioassay. *Proc. Roy. Soc. London Ser. B* **145** 42–61.
[5] Grizzle, J. E. (1965). The two-period change-over design and its use in clinical trials. *Biometrics* **21** 467–480.
[6] Hedayat, A. S. and Afsarinejad, K. (1978). Repeated measurements designs. II. *Ann. Statist.* **6** 619–628. MR0488527
[7] Hedayat, A. S. and Yang, M. (2003). Universal optimality of balanced uniform crossover designs. *Ann. Statist.* **31** 978–983. MR1994737
[8] Keppel, G. (1973). *Design and Analysis*: *A Researcher's Handbook*. Prentice–Hall, Englewood Cliffs, NJ.
[9] Kiefer, J. (1975). Construction and optimality of generalized Youden designs. In *A Survey of Statistical Design and Linear Models* (J. N. Srivastava, ed.) 333–353. North-Holland, Amsterdam. MR0395079
[10] Kunert, J. (1984). Optimality of balanced uniform repeated measurements designs. *Ann. Statist.* **12** 1006–1017. MR0751288
[11] Kunert, J. and Martin, R. J. (2000). On the determination of optimal designs for an interference model. *Ann. Statist.* **28** 1728–1741. MR1835039

[12] KUNERT, J. and MARTIN, R. J. (2000). Optimality of type I orthogonal arrays for crossover models with correlated errors. *J. Statist. Plann. Inference* **87** 119–124. [MR1772043](MR1772043)

[13] KUNERT, J. and STUFKEN, J. (2002). Optimal crossover designs in a model with self and mixed carryover effects. *J. Amer. Statist. Assoc.* **97** 898–906. [MR1941418](MR1941418)

[14] KUSHNER, H. B. (1997). Optimal repeated measurements designs: The linear optimality equations. *Ann. Statist.* **25** 2328–2344. [MR1604457](MR1604457)

[15] KUSHNER, H. B. (1998). Optimal and efficient repeated-measurements designs for uncorrelated observations. *J. Amer. Statist. Assoc.* **93** 1176–1187. [MR1649211](MR1649211)

[16] MARKIEWICZ, A. (1997). Properties of information matrices for linear models and universal optimality of experimental designs. *J. Statist. Plann. Inference* **59** 127–137. [MR1450780](MR1450780)

[17] MIELKE, P. W., JR. (1974). Square rank test appropriate to weather modification cross-over design. *Technometrics* **16** 13–16. [MR0343498](MR0343498)

[18] PATTERSON, H. D. (1952). The construction of balanced designs for experiments involving sequences of treatments. *Biometrika* **39** 32–48. [MR0050855](MR0050855)

[19] PUKELSHEIM, F. (1993). *Optimal Design of Experiments*. Wiley, New York. [MR1211416](MR1211416)

[20] SHAH, K. R. and SINHA, B. K. (1989). *Theory of Optimal Designs. Lecture Notes in Statist.* **54**. Springer, New York. [MR1016151](MR1016151)

[21] SHAH, K. R. and SINHA, B. K. (2002). Universal optimality for the joint estimation of parameters. Unpublished manuscript.

[22] STUFKEN, J. (1996). Optimal crossover designs. In *Design and Analysis of Experiments* (S. Ghosh and C. R. Rao, eds.) 63–90. North-Holland, Amsterdam. [MR1492565](MR1492565)

[23] WILLIAMS, E. J. (1949). Experimental designs balanced for the estimation of residual effects of treatments. *Australian J. Sci. Research. Ser. A* **2** 149–168. [MR0033508](MR0033508)

K. R. SHAH
DEPARTMENT OF STATISTICS
AND ACTUARIAL SCIENCE
UNIVERSITY OF WATERLOO
WATERLOO, ONTARIO
CANADA N2L 361
E-MAIL: [dakshakirti@yahoo.ca](dakshakirti@yahoo.ca)

M. BOSE
APPLIED STATISTICS UNIT
INDIAN STATISTICAL INSTITUTE
203 B. T. ROAD
KOLKATA 7000108
INDIA
E-MAIL: [mausumi@isical.ac.in](mausumi@isical.ac.in)

D. RAGHAVARAO
DEPARTMENT OF STATISTICS
TEMPLE UNIVERSITY
PHILADELPHIA, PENNSYLVANIA 19122-6083
USA
E-MAIL: [draghava@temple.edu](draghava@temple.edu)